\pdfoutput=1
\RequirePackage{ifpdf}
\ifpdf 
\documentclass[pdftex]{sigma}
\else
\documentclass{sigma}
\fi

\usepackage[all]{xy}
\numberwithin{equation}{section}
\newtheorem*{theointro}{Theorem}
\newtheorem{Theorem}{Theorem}[section]
\newtheorem{Corollary}[Theorem]{Corollary}
\newtheorem{Lemma}[Theorem]{Lemma}
\newtheorem{Proposition}[Theorem]{Proposition}
{ \theoremstyle{definition}
\newtheorem{Definition}[Theorem]{Definition}
\newtheorem{Remark}[Theorem]{Remark} }

\begin{document}

\newcommand{\divi}[1]{\operatorname{div}_{#1}}
\newcommand{\weight}[1]{\omega_{#1}}
\newcommand{\degr}[1]{\deg_{#1}}
\newcommand{\wpf}[1]{\wp_{#1}}
\newcommand{\prodtheta}[1]{\Theta_{#1}}
\newcommand{\quotprodtheta}[1]{\Theta_{#1}^{\rm quot}}

\newcommand{\arXivNumber}{1405.2002}

\allowdisplaybreaks

\renewcommand{\thefootnote}{$\star$}

\renewcommand{\PaperNumber}{003}

\FirstPageHeading

\ShortArticleName{Galois Groups of Dif\/ference Equations of Order Two on Elliptic Curves}

\ArticleName{Galois Groups of Dif\/ference Equations\\
of Order Two on Elliptic Curves\footnote{This paper is a~contribution to the Special Issue on Algebraic Methods in
Dynamical Systems. The full collection is available at
\href{http://www.emis.de/journals/SIGMA/AMDS2014.html}{http://www.emis.de/journals/SIGMA/AMDS2014.html}}}

\Author{Thomas DREYFUS~$^\dag$ and Julien ROQUES~$^\ddag$}

\AuthorNameForHeading{T.~Dreyfus and J.~Roques}

\Address{$^\dag$~Universit\'e Paul Sabatier - Institut de Math\'ematiques de Toulouse,\\
\hphantom{$^\dag$}~18 route de Narbonne, 31062 Toulouse, France}
\EmailD{\href{mailto:tdreyfus@math.univ-toulouse.fr}{tdreyfus@math.univ-toulouse.fr}}
\URLaddressD{\url{https://sites.google.com/site/thomasdreyfusmaths/}}

\Address{$^\ddag$~Institut Fourier, Universit\'e Grenoble 1, CNRS UMR 5582, 100 rue des Maths,\\
\hphantom{$^\ddag$}~BP 74, 38402 St Martin d'H\`eres, France}
\EmailD{\href{mailto:Julien.Roques@ujf-grenoble.fr}{Julien.Roques@ujf-grenoble.fr}}
\URLaddressD{\url{www-fourier.ujf-grenoble.fr/~jroques/}}

\ArticleDates{Received August 06, 2014, in f\/inal form January 08, 2015; Published online January 13, 2015}

\Abstract{This paper is concerned with dif\/ference equations on elliptic curves.
We establish some general properties of the dif\/ference Galois groups of equations of order two, and give applications to
the calculation of some dif\/ference Galois groups.
For instance, our results combined with a~result from transcendence theory due to Schneider allow us to identify a~large
class of discrete Lam\'e equations with dif\/ference Galois group $\operatorname{GL}_{2}(\mathbb C)$.}

\Classification{39A06; 12H10}

\renewcommand{\thefootnote}{\arabic{footnote}}
\setcounter{footnote}{0}

\section{Introduction}

Let $\mathcal{E} \subset \mathbb{P}^{2}$ be the elliptic curve def\/ined by the projectivization of the Weierstrass
equation
\begin{gather}
\label{Weierstrass equation}
y^{2}=4x^{3}-g_{2}x-g_{3} \qquad \hbox{with} \quad g_{2},g_{3} \in \mathbb C.
\end{gather}
We denote by $\mathcal{E}(\mathbb C)$ the group of $\mathbb C$-points of $\mathcal{E}$.
Its abelian group law is denoted by $\oplus$.

In this paper, we study the dif\/ference Galois groups of linear dif\/ference equations of order two on $\mathcal{E}(\mathbb
C)$ of the form:
\begin{gather}
\label{equa sur courbe elliptique generique}
\underline{y}(\underline{z} \oplus 2\underline{h}) + \underline{a}(\underline{z}) \underline{y}(\underline{z}\oplus
\underline{h}) + \underline{b}(\underline{z}) \underline{y}(\underline{z}) =0,
\end{gather}
where $\underline{y}$ is an unknown function of the variable $\underline{z} \in \mathcal{E}(\mathbb C)$, $\underline{h}$
is a~f\/ixed non torsion point of $\mathcal{E}(\mathbb C)$ and $\underline{a}$, $\underline{b}$ are given rational functions
on~$\mathcal{E}$.

This equation can be seen as a~dif\/ference equation over $\mathbb C$.
Indeed, if $\Lambda \subset \mathbb C$ is a~lattice of periods of $\mathcal{E}$ and if $\wp$ is the corresponding
Weierstrass function, then $\mathbb C/\Lambda$ is identif\/ied with $\mathcal{E}(\mathbb C)$ via the factorization through
$\mathbb C/\Lambda$ of
\begin{gather*}
\varphi: \ z \in \mathbb C \mapsto (\wp(z):\wp'(z):1) \in \mathcal{E}(\mathbb C).
\end{gather*}
Pulling back the equation~\eqref{equa sur courbe elliptique generique} via~$\varphi$, which is a~group morphism from
$(\mathbb C, +)$ to $(\mathcal{E}(\mathbb C),\oplus)$, we obtain the following dif\/ference equation on $\mathbb C$:
\begin{gather}
\label{equa intro}
y(z+2h) + a(z) y(z+h) + b(z)y(z) =0,
\end{gather}
where~$y$ is an unknown function of the variable $z\in \mathbb C$, $a:=\underline{a} \circ \varphi$ and
$b:=\underline{b} \circ \varphi$ are~$\Lambda$-periodic elliptic functions and ${h \in \varphi^{-1}(\underline{h})}$.
So, we have the following relations between the equations~\eqref{equa sur courbe elliptique generique} and~\eqref{equa
intro}: $\underline{z} = \varphi(z)$, $\underline{h}=\varphi(h)$ and $\underline{y}(\underline{z})= y(z)$.

These equations are discrete counterparts of dif\/ferential equations on elliptic curves, a~famous example of which is
Lam\'e dif\/ferential equation
\begin{gather*}
y''(z) = (A\wp(z)+B)y(z),
\end{gather*}
where $A,B\in \mathbb C$.
The main results of this paper allow us to compute the dif\/ference Galois groups of some equations such as the discrete
Lam\'e equation
\begin{gather}
\label{equa lame deform intro}
\Delta_h^2 y = (A\wp(z)+B)y,\qquad \text{where} \quad \Delta_h y(z) = \frac{y(z+h)-y(z)}{h}.
\end{gather}
For instance, the following theorem is a~consequence of our main results combined with a~result from transcendence
theory due to Schneider in~\cite{schneider} (see also Bertrand and Masser's pa\-pers~\mbox{\cite{BerMas,Mas75}}).

\begin{theointro}
Assume that $\mathcal{E}$ is defined over $\overline{\mathbb{Q}}$ $($i.e., $g_{2},g_{3} \in \overline{\mathbb{Q}})$
and that $h,A,B\in \overline{\mathbb{Q}}$ with $A \neq 0$.
Then, the difference Galois group of equation~\eqref{equa lame deform intro} is $\operatorname{GL}_{2}(\mathbb C)$.
\end{theointro}

To be precise, the base f\/ield for the dif\/ference Galois groups considered in the present paper is not the f\/ield
of~$\Lambda$-periodic meromorphic functions over $\mathbb C$, but the f\/ield constituted of the meromorphic functions
over $\mathbb C$ which are $\Lambda'$-periodic for some sub-lattice $\Lambda'$ of~$\Lambda$.

The galoisian aspects of the theory of dif\/ference equations have attracted the attention of many authors in the past
years,
e.g.,~\cite{andrenoncomm,BuT,chatziharsing,diviziogrothendieck,hardouin-divizioqgrothendieckCRAS,etingofgalois,franke63,franke66,franke67,franke74,
hardouin-singer,NG,OW14,ramissauloyqstokes1,ramissauloyqstokes2,RS14,sauloyqgaloisfuchs,vdpr,VdPS97}.
The calculation of the dif\/ference Galois groups of f\/inite dif\/ference or~$q$-dif\/ference equations of order two on
$\mathbb{P}^{1}$ has been considered by Hendricks~\cite{Hen97,Hen98} and by the second author~\cite{Ro08}.
The work of Hendricks served as a~basis for the present work, but, to the best of our knowledge, the present paper is
the f\/irst to consider the dif\/ference Galois groups of dif\/ference equations on a~non rational variety.
The study of dynamical systems on elliptic curves appears in several areas of mathematics (e.g., discrete dynamical
systems, QRT maps).
In particular, it is very likely that the equations considered in the present paper will arise as linearizations of
discrete dynamical systems, in connection with discrete Morales--Ramis theories~\cite{CasR,casale-roques-2}.
In this context, the dif\/ference Galois groups are used to obtain non-integrability results.

This paper is organized as follows.
Section~\ref{sec2} contains reminders and complements on dif\/ference Galois theory (for equations of arbitrary order)
with a~special emphasis on dif\/ference equations on elliptic curves.
We insist on the fact that the base dif\/ference f\/ield for the dif\/ference Galois groups considered in the present paper is
not the f\/ield of~$\Lambda$-periodic elliptic functions but the f\/ield of elliptic functions which are $\Lambda'$-periodic
for some sub-lattice $\Lambda'$ of~$\Lambda$.
In Section~\ref{sec3}, we introduce some notations related to the special functions used in this paper (theta
functions, Weierstrass $\wp$-functions) and we collect some useful results.
In Section~\ref{sec4}, we study the relations between the irreducibility of the dif\/ference Galois group of
equation~\eqref{equa intro} and the solutions of an associated Riccati-type equation.
We then study this Riccati equation assuming that we have {\it a~priori} informations on the divisors of the
coef\/f\/icients~$a$ and~$b$.
In Section~\ref{sec5}, we show that there is a~similar relation between the imprimitivity of the Galois group and some
Riccati-type equation.
Section~\ref{sec6} is devoted to the calculation of some dif\/ference Galois groups, including those of the discrete
Lam\'e equations mentioned above.

\section{Dif\/ference Galois theory: reminders and complements}
\label{sec2}

\subsection{Generalities on dif\/ference Galois theory}

For details on what follows, we refer to~\cite[Chapter~1]{VdPS97}.

A dif\/ference ring $(R,\phi)$ is a~ring~$R$ together with a~ring automorphism $\phi: R \rightarrow R$.
An ideal of~$R$ stabilized by~$\phi$ is called a~dif\/ference ideal of $(R,\phi)$.
If~$R$ is a~f\/ield then $(R,\phi)$ is called a~dif\/ference f\/ield.

The ring of constants $R^{\phi}$ of the dif\/ference ring $(R,\phi)$ is def\/ined~by
\begin{gather*}
R^{\phi}:=\{f \in R \, | \, \phi(f)=f\}.
\end{gather*}

A dif\/ference ring morphism (resp.\
dif\/ference ring isomorphism) from the dif\/ference ring $(R,\phi)$ to the dif\/ference ring
$(\widetilde{R},\widetilde{\phi})$ is a~ring morphism (resp.\
ring isomorphism) $\varphi: R \rightarrow \widetilde{R}$ such that $\varphi \circ \phi = \widetilde{\phi} \circ
\varphi$.

A dif\/ference ring $(\widetilde{R},\widetilde{\phi})$ is a~dif\/ference ring extension of a~dif\/ference ring $(R,\phi)$ if
$\widetilde{R}$ is a~ring extension of~$R$ and $\widetilde{\phi}_{\vert R}=\phi$; in this case, we will often denote
$\widetilde{\phi}$ by~$\phi$.
Two dif\/ference ring extensions $(\widetilde{R}_{1},\widetilde{\phi}_{1})$ and $(\widetilde{R}_{2},\widetilde{\phi}_{2})$
of a~dif\/ference ring $(R,\phi)$ are isomorphic over $(R,\phi)$ if there exists a~dif\/ference ring isomorphism~$\varphi$
from $(\widetilde{R}_{1},\widetilde{\phi}_{1})$ to $(\widetilde{R}_{2},\widetilde{\phi}_{2})$ such that $\varphi_{\vert
R}=\operatorname{Id}_{R}$.

{\it We now let $(K,\phi)$ be a~difference field.
We assume that its field of constants $C:=K^{\phi}$ is algebraically closed and that the characteristic of~$K$ is~$0$.}

Consider a~linear dif\/ference system
\begin{gather}
\label{the generic system}
\phi Y=AY \qquad \hbox{with} \quad A \in \operatorname{GL}_{n}(K).
\end{gather}

According to~\cite[Section~1.1]{VdPS97}, there exists a~dif\/ference ring extension $(R,\phi)$ of $(K,\phi)$ such that
\begin{itemize}\itemsep=0pt
\item[1)] there exists $U \in \operatorname{GL}_{n}(R)$ such that $\phi (U) = AU$ (such a~$U$ is called a~fundamental
matrix of solutions of~\eqref{the generic system});
\item[2)] $R$ is generated, as a~$K$-algebra, by the entries of~$U$ and $\det(U)^{-1}$;
\item[3)] the only dif\/ference ideals of $(R,\phi)$ are $\{0\}$ and~$R$.
\end{itemize}
Such a~dif\/ference ring $(R,\phi)$ is called a~Picard--Vessiot ring for~\eqref{the generic system} over~$(K,\phi)$.
It is unique up to isomorphism of dif\/ference rings over $(K,\phi)$.
It is worth mentioning that $R^{\phi}=C$; see~\cite[Lemma~1.8]{VdPS97}.
\begin{Remark}
Picard--Vessiot rings are not domains in general: they are f\/inite direct sums of domains cyclically permuted by~$\phi$;
see~\cite[Corollary~1.16]{VdPS97}.
\end{Remark}

The corresponding dif\/ference Galois group~$G$ over $(K,\phi)$ of~\eqref{the generic system} is the group of~$K$-linear
ring automorphisms of~$R$ commuting with~$\phi$:
\begin{gather*}
G:=\{\sigma \in \operatorname{Aut}(R/K) \, | \, \phi\circ\sigma=\sigma\circ \phi \}.
\end{gather*}
The choice of the base f\/ield is by no way innocent.
The bigger the base f\/ield is, the smaller the Galois group is.

A straightforward computation shows that, for any~$\sigma {\in} G$, there exists a~unique ${C(\sigma) {\in}
\operatorname{GL}_{n}(C)}$ such that $\sigma(U)=UC(\sigma)$.
According to~\cite[Theorem~1.13]{VdPS97}, one can identify~$G$ with an {\it algebraic} subgroup of
$\operatorname{GL}_{n}(C)$ via the faithful representation
\begin{gather*}
\sigma \in G \mapsto C(\sigma) \in \operatorname{GL}_{n}(C).
\end{gather*}
If we choose another fundamental matrix of solutions~$U$, we f\/ind a~conjugate representation.

\begin{Remark}
Given an~$n$th order dif\/ference equation
\begin{gather}
\label{the generic equa}
a_{n} \phi^{n} (y) + \cdots + a_{1} \phi(y)+a_{0} y=0,
\end{gather}
with $a_{0},\dots,a_{n} \in K$ and $a_{0}a_{n} \neq 0$, we can consider the equivalent linear dif\/ference system
\begin{gather}
\label{matrix system equa}
\phi Y=AY, \qquad \text{with} \quad  A=
\begin{pmatrix}
0&1&0&\dots&0
\\
0&0&1&\ddots&\vdots
\\
\vdots&\vdots&\ddots&\ddots&0
\\
0&0&\dots&0&1
\\
-\frac{a_{0}}{a_{n}}& -\frac{a_{1}}{a_{n}}&\dots & \dots & -\frac{a_{n-1}}{a_{n}}
\end{pmatrix}
\in \operatorname{GL}_{n}(K).
\end{gather}
By ``Galois group of the dif\/ference equation~\eqref{the generic equa}'' we mean ``Galois group of the dif\/ference
system~\eqref{matrix system equa}''.
\end{Remark}

We shall now introduce a~property relative to the dif\/ference base f\/ield, which appears in~\cite[Lemma~1.19]{VdPS97}.

\begin{Definition}
\label{prop P}
We say that the dif\/ference f\/ield $(K,\phi)$ satisf\/ies property $(\mathcal{P})$ if the following properties hold:
\begin{itemize}\itemsep=0pt
\item the f\/ield~$K$ is a~$\mathcal{C}^{1}$-f\/ield\footnote{Recall that~$K$ is a~$\mathcal{C}^{1}$-f\/ield if every
non-constant homogeneous polynomial~$P$ over~$K$ has a~non-trivial zero provided that the number of its variables is
more than its degree.};
\item if~$L$ is a~f\/inite f\/ield extension of~$K$ such that~$\phi$ extends to a~f\/ield endomorphism of~$L$ then $L=K$.
\end{itemize}
\end{Definition}

The following result is due to van der Put and Singer.
We recall that two dif\/ference systems $\phi Y=AY$ and $\phi Y=BY$ with $A,B \in \operatorname{GL}_{n}(K)$ are isomorphic
over~$K$ if and only if there exists $T \in \operatorname{GL}_{n}(K)$ such that $\phi(T) A=BT$.

\begin{Theorem}
\label{si prop P alors}
Assume that $(K,\phi)$ satisfies property $(\mathcal{P})$.
Let $K^{\phi} = C$.
Let $G \subset \operatorname{GL}_{n}(C)$ be the difference Galois group over $(K,\phi)$ of
\begin{gather}
\label{the generic system theo}
\phi(Y) = AY, \qquad \text{with} \quad  A \in \operatorname{GL}_{n}(K).
\end{gather}
Then, the following properties hold:
\begin{itemize}
\itemsep=0pt
\item $G/G^{\circ}$ is cyclic, where $G^{\circ}$ is the identity component of~$G$;
\item there exists $B \in G(K)$ such that~\eqref{the generic system theo} is isomorphic to $\phi Y=BY$ over~$K$.
\end{itemize}
Let $\widetilde{G}$ be an algebraic subgroup of $\operatorname{GL}_{n}(C)$ such that~$A\in \widetilde{G}(K)$.
The following properties hold:
\begin{itemize}
\itemsep=0pt
\item~$G$ is conjugate to a~subgroup of $\widetilde{G}$;
\item any minimal element in the set of algebraic subgroups $\widetilde{H}$ of $\widetilde{G}$ for which there exists
$T\in \operatorname{GL}_{n}(K)$ such that $\phi(T)AT^{-1}\in \widetilde{H}(K)$ is conjugate to~$G$;
\item~$G$ is conjugate to $\widetilde{G}$ if and only if, for any $T\in \widetilde{G}(K)$ and for any proper algebraic
sub\-group~$\widetilde{H}$ of~$\widetilde{G}$, one has that $\phi(T)AT^{-1}\notin \widetilde{H}(K)$.
\end{itemize}
\end{Theorem}

\begin{proof}
The proof of~\cite[Propositions~1.20 and~1.21]{VdPS97} in the special case where $K:=C(z)$ and~$\phi$ is the shift
$\phi(f(z)):=f(z+h)$ with $h \in C^\times$, extends {\it mutatis mutandis} to the present case.
\end{proof}

\subsection{Dif\/ference equations on elliptic curves}

Let $\Lambda \subset \mathbb C$ be a~lattice.
Without loss of generality, we can assume that
\begin{gather*}
\Lambda=\mathbb Z +\mathbb Z \tau, \qquad \text{with} \quad \Im (\tau)>0,
\end{gather*}
where $\Im(\cdot)$ denotes the imaginary part.
For any lattice $\Lambda' \subset \mathbb C$, we let $\operatorname{M}_{\Lambda'}$ be the f\/ield of $\Lambda'$-periodic
meromorphic functions.
We denote by~$K$ the f\/ield def\/ined by
\begin{gather*}
K:= \bigcup_{\Lambda' \; \text{sub-lattice of} \; \Lambda} \operatorname{M}_{\Lambda'} =\bigcup_{k \geq 1} \operatorname{M}_{k  \Lambda}.
\end{gather*}
Let $h \in \mathbb C$ such that $h \mod \Lambda$ is not a~torsion point of $\mathbb C/\Lambda$.
We endow~$K$ with the non-cyclic f\/ield automorphism~$\phi$ def\/ined~by
\begin{gather*}
\phi(f)(z):=f(z+h).
\end{gather*}
Then,~$(K,\phi)$ is a~dif\/ference f\/ield.

\begin{Proposition}
The field of constants of $(K,\phi)$ is
\begin{gather*}
K^{\phi}=\mathbb C.
\end{gather*}
\end{Proposition}

\begin{proof}
Consider $f\in K^{\phi}$.
Let $\Lambda'$ be a~sub-lattice of~$\Lambda$ such that $f \in \operatorname{M}_{\Lambda'}$.
Note that~$f$ is $\Lambda'$-periodic (because $f \in \operatorname{M}_{\Lambda'}$) and~$h$-periodic (because
$\phi(f)=f$), so~$f$ is a~$(\Lambda' + h\mathbb Z)$-periodic meromorphic function.
But $\Lambda' + h\mathbb Z$ has an accumulation point because $h \mod \Lambda$ is not a~torsion point of $\mathbb
C/\Lambda$.
Therefore,~$f$ is constant.
\end{proof}

\begin{Proposition}
The difference field $(K,\phi)$ satisfies property $(\mathcal{P})$ $($see Definition~{\rm \ref{prop P})}.
\end{Proposition}

\begin{proof}
Since $K = \bigcup_{k \geq 1} \operatorname{M}_{k \Lambda}$ is the increasing union of the f\/ields
$\operatorname{M}_{k \Lambda}$, the fact that~$K$ is a~$\mathcal{C}^{1}$-f\/ield follows from Tsen's
theorem~\cite{LangQAC} (according to which the function f\/ield of any algebraic curve over an algebraically closed f\/ield,
e.g.,
$\operatorname{M}_{k \Lambda}$, is $\mathcal{C}^{1}$).

Let~$L$ be a~f\/inite extension of~$K$ such that~$\phi$ extends to a~f\/ield endomorphism of~$L$.
We have to prove that $L=K$.
The primitive element theorem ensures that there exists $u \in L$ such that $L=K(u)$.
Let $\Lambda'$ be a~sub-lattice of~$\Lambda$ such that
\begin{itemize}
\itemsep=0pt
\item~$u$ is algebraic over $\operatorname{M}_{\Lambda'}$,
\item $\phi(u) \in \operatorname{M}_{\Lambda'}(u)$.
\end{itemize}
Then, $\operatorname{M}_{\Lambda'}(u)$ is a~f\/inite extension of $\operatorname{M}_{\Lambda'}$ and~$\phi$ induces an
automorphism of $\operatorname{M}_{\Lambda'}(u)$.

Using the equivalence of categories between between smooth projective curves and function f\/ields of dimension
$1$~\cite[Corollary 6.12]{HartshorneAG}, we see that there exists a~commutative diagram of the form
\begin{gather*}
\xymatrix{X \ar[r]^{f} \ar[d]_{\varphi} &  X \ar[d]^{\varphi}
\\
\mathbb C/\Lambda' \ar[r]_{z \mapsto z+h}  & \mathbb C/\Lambda'}
\end{gather*}
where $\varphi: X \rightarrow \mathbb C/\Lambda'$ is a~morphism of smooth projective curves, whose induced morphism of
function f\/ields ``is'' the inclusion $\operatorname{M}_{\Lambda'} \subset \operatorname{M}_{\Lambda'}(u)$, and where~$f$
is an endomorphism of~$X$, whose induced morphism on function f\/ields ``is'' $\phi: \operatorname{M}_{\Lambda'}(u)
\rightarrow \operatorname{M}_{\Lambda'}(u)$.
Considering this commutative diagram, we see that~$f$ has degree $1$ and that, if~$\varphi$ is ramif\/ied above $y \in
\mathbb C/\Lambda'$, then~$\varphi$ is also ramif\/ied above $y-h$.
So, the set of ramif\/ication values of~$\varphi$ is stable by $z \mapsto z-h$.
This set being f\/inite, it has to be empty.
So,~$\varphi$ is unramif\/ied.

Hurwitz's formula implies that~$X$ has genus~$1$, i.e., that~$X$ is an elliptic curve.
So, there exist a~lattice $\Lambda'' \subset \mathbb C$ and an isomorphism $\psi: \operatorname{M}_{\Lambda'}(u)
\rightarrow \operatorname{M}_{\Lambda''}$.
There exists $(a,b) \in \mathbb C^{\times} \times \mathbb C$ such that $a \Lambda'' \subset \Lambda'$ and such that the
restriction $\psi_{\vert \operatorname{M}_{\Lambda'}}$ is given, for all $f \in \operatorname{M}_{\Lambda'}$,~by
$\psi(f)(z)=f(az+b)$.
(Indeed, $\psi_{\vert \operatorname{M}_{\Lambda'}}: \operatorname{M}_{\Lambda'} \rightarrow
\operatorname{M}_{\Lambda''}$ is a~f\/ield morphism from the function f\/ield of the elliptic curve~$\mathbb C/\Lambda'$ to
the function f\/ield of the elliptic curve~$\mathbb C/\Lambda''$.
So, $\psi_{\vert \operatorname{M}_{\Lambda'}}$ is induced by a~morphism from the elliptic curve $\mathbb C/\Lambda''$ to
the elliptic curve $\mathbb C/\Lambda'$.
Now, our claim follows from the fact that the morphisms from $\mathbb C/\Lambda''$ to $\mathbb C/\Lambda'$ are of the
form $z \mod \Lambda'' \mapsto az+b \mod \Lambda'$ for some $(a,b) \in \mathbb C^{\times} \times \mathbb C$ such that $a
\Lambda'' \subset \Lambda'$.) The commutative diagram
\begin{gather*}
\xymatrix{\operatorname{M}_{\Lambda'} \ar@{^{(}->}[r] \ar[rd]_{v(z) \mapsto v(az+b)} &  \operatorname{M}_{\Lambda'}(u)
\ar[d]^{\psi} \ar[rd] &
\\
&  \operatorname{M}_{\Lambda''} \ar[r]_{w(z) \mapsto w(\frac{z-b}{a})}  & \operatorname{M}_{a \Lambda''}}
\end{gather*}
shows that the f\/ields $\operatorname{M}_{\Lambda'}(u)$ and $\operatorname{M}_{a \Lambda''}$ are
$\operatorname{M}_{\Lambda'}$-isomorphic.
But the extension $\operatorname{M}_{a \Lambda''}/\operatorname{M}_{\Lambda'}$ is Galois (indeed, this is equivalent to
the fact that the corresponding morphism of smooth projective curves $\mathbb C/\Lambda' \rightarrow \mathbb
C/a\Lambda''$ is Galois, and this is easily seen from the explicit description of the morphisms between these curves).
Therefore, any $\operatorname{M}_{\Lambda'}$-morphism from $\operatorname{M}_{a \Lambda''}$ to $K (u)$ must leave
$\operatorname{M}_{a \Lambda''}$ globally invariant.
But, $\operatorname{M}_{a \Lambda''}$ and $\operatorname{M}_{\Lambda'}(u)$ are $\operatorname{M}_{\Lambda'}$-isomorphic
subf\/ields of~$K(u)$.
So $\operatorname{M}_{\Lambda'} (u) \subset \operatorname{M}_{a \Lambda''}$, and therefore $u \in \operatorname{M}_{a
\Lambda''} \subset K$ and $L=K(u) \subset K$.
\end{proof}

\begin{Corollary}
\label{coro1}
The conclusions of Theorem~{\rm \ref{si prop P alors}} are valid for $(K,\phi)$.
\end{Corollary}

\section[Theta functions and Weierstrass $\wp$-function]{Theta functions and Weierstrass $\boldsymbol{\wp}$-function}
\label{sec3}

\subsection{Theta functions}
\label{sec31}

We recall that
\begin{gather*}
\Lambda=\mathbb Z+\tau \mathbb Z \subset \mathbb C \qquad \text{with} \quad \Im(\tau)>0.
\end{gather*}

Let~$\theta$ be the Jacobi theta function def\/ined~by
\begin{gather*}
\theta(z)=\sum\limits_{m\in \mathbb Z}(-1)^{m} e^{i\pi m(m-1) \tau} e^{2i \pi m z}.
\end{gather*}
We shall now recall some basic facts about this function.
We refer to~\cite[Chapter~I]{Mumford} for details and proofs.

\begin{Remark}
The classical theta function is def\/ined by ${\vartheta(z,\tau) = \sum\limits_{m \in \mathbb Z} e^{i\pi m^{2} \tau+ 2i
\pi m z}}$.
Actually, this is the function studied in~\cite[Chapter~I]{Mumford}.
But, there is a~simple relation between~$\theta$ and~$\vartheta$, namely $\theta(z)=\vartheta(z+\frac{1-\tau}{2},\tau)$.
So that any statement for~$\vartheta$ can be immediately translated into a~statement for~$\theta$.
\end{Remark}

We recall that~$\theta$ is a~$1$-periodic entire function such that
\begin{gather*}
\theta(z+\tau)=- e^{-2i\pi z}\theta(z).
\end{gather*}
Moreover, we have the following formula, known as Jacobi's triple product formula:
\begin{gather*}
\theta(z)= \prod\limits_{m=1}^{\infty}\big(1-e^{2i\pi \tau
m}\big)\big(1-e^{2i\pi((m-1)\tau+z)}\big)\big(1-e^{2i\pi(m\tau-z)}\big).
\end{gather*}
For any integer $k \geq 1$, we let $\theta_{k}$ be the function given~by
\begin{gather*}
\theta_{k}(z):=\theta(z/k).
\end{gather*}
This~$k$-periodic entire function satisf\/ies the following functional equation:
\begin{gather}
\label{func equa theta k}
\theta_{k}(z+k\tau)=- e^{-2i\pi z/k}\theta_{k}(z).
\end{gather}
It follows from Jacobi's triple product formula that the zeroes of $\theta_{k}$ are simple and that its set of zeroes is
$k\Lambda$.

Let $\prodtheta{k}$ be the set of entire functions of the form
\begin{gather*}
c  \prod\limits_{\xi \in \mathbb C} \theta_{k}(z-\xi)^{n_{\xi}}
\end{gather*}
with $c \in \mathbb C^{\times}$ and $(n_{\xi})_{\xi \in \mathbb C}\in \mathbb N^{(\mathbb C)}$ with f\/inite support.
We denote by $\quotprodtheta{k}$ the set of meromorphic functions over $\mathbb C$ that can be written as quotient of
two elements of $\prodtheta{k}$.

We def\/ine the divisor $\divi{k} (f)$ of $ f \in \quotprodtheta{k} $ as the following formal sum of points of $\mathbb
C/k\Lambda$:
\begin{gather*}
\divi{k} (f):= \sum\limits_{\lambda \in \mathbb C/k\Lambda} \operatorname{ord}_{\lambda}(f) [\lambda],
\end{gather*}
where $\operatorname{ord}_{\lambda}(f)$ is the $(z-\xi)$-adic valuation of~$f$, for an arbitrary $\xi \in \lambda$ (it
does not depend on the chosen $\xi \in \lambda$).
For any $\lambda \in \mathbb C/k\Lambda$ and any~$\xi \in \lambda$, we set
\begin{gather*}
[\xi]_k:=[\lambda].
\end{gather*}
Moreover, we will write
\begin{gather*}
\sum\limits_{\lambda \in \mathbb C/k\Lambda} n_{\lambda} [\lambda] \leq \sum\limits_{\lambda \in \mathbb C/k\Lambda}
m_{\lambda} [\lambda]
\end{gather*}
if, for all $\lambda \in \mathbb C/k\Lambda$, $n_{\lambda} \leq m_{\lambda}$.
We also introduce the weight $\weight{k} (f)$ of~$f$ def\/ined~by
\begin{gather*}
\weight{k}(f):=\sum\limits_{\lambda \in \mathbb C/k\Lambda} \operatorname{ord}_{\lambda}(f) \lambda \in \mathbb
C/k\Lambda
\end{gather*}
and its degree $\degr{k} (f)$ given~by
\begin{gather*}
\degr{k} (f):=\sum\limits_{\lambda \in \mathbb C/k\Lambda} \operatorname{ord}_{\lambda}(f) \in \mathbb Z.
\end{gather*}

If $f= c  \prod\limits_{\xi \in \mathbb C} \theta_{k}(z-\xi)^{n_{\xi}} \in \quotprodtheta{k}$, then
\begin{gather*}
\divi{k} (f) = \sum\limits_{\xi \in \mathbb C} n_{\xi} [\xi]_{k},
\qquad
\weight{k}(f) = \sum\limits_{\xi \in \mathbb C} n_\xi \xi \mod k\Lambda
\qquad {\rm and} \qquad
\degr{k} (f)=\sum\limits_{\xi \in \mathbb C} n_\xi.
\end{gather*}

The interest of $\quotprodtheta{k}$ in our context is given by the following classical result.

\begin{Proposition}
We have
\begin{gather*}
\operatorname{M}_{k\Lambda}^{\times} \subset \quotprodtheta{k}.
\end{gather*}
\end{Proposition}

\begin{proof}
This inclusion means that any $k\Lambda$-periodic meromorphic function can be written, up to some multiplicative
constant in $\mathbb C^{\times}$, as a~quotient of product of functions of the form $\theta_{k}(z-\xi)$.
This is classical, see~\cite[Chapter I, Section~6]{Mumford}.
\end{proof}

We now state a~couple of lemmas, which will be used freely in the rest of the paper.

\begin{Lemma}
\label{form precise}
Any $f= c  \prod\limits_{\xi \in \mathbb C} \theta_{k}(z-\xi)^{n_{\xi}} \in \quotprodtheta{k}$
is~$k$-periodic and satisfies
\begin{gather}
\label{caract theta quot}
f(z+k\tau)=(-1)^{\degr{k} (f)} e^{2i\pi \omega} e^{-2i\pi \degr{k} (f) z/k}f(z),
\end{gather}
where   $\omega = \sum\limits_{\xi \in \mathbb C} n_\xi \xi$ is a~representative of $\weight{k}(f)$\footnote{It follows
from this formula that~$f$ belongs to $\operatorname{M}_{k\Lambda}$ if and only if $\degr{k} (f)=\sum\limits_{\xi \in
\mathbb C} n_\xi=0$ and $\omega=\sum\limits_{\xi \in \mathbb C} n_\xi \xi \in \mathbb Z$.}.
Conversely, any non zero~$k$-periodic meromorphic function~$f$ over $\mathbb C$ such that
\begin{gather}
\label{caract theta quot bis}
f(z+k\tau)=c e^{-2i\pi n z/k}f(z),
\end{gather}
for some $c \in \mathbb C^{\times}$ and $n \in \mathbb Z$, belongs to~$\quotprodtheta{k}$.
\end{Lemma}

\begin{proof}
The fact that any $f \in \quotprodtheta{k}$ is~$k$-periodic and satisf\/ies the functional equation~\eqref{caract theta
quot} follows from the fact that $\theta_{k}$ is~$k$-periodic and satisf\/ies the functional equation~\eqref{func equa
theta k}.
Conversely, consider a~non zero~$k$-periodic meromorphic function~$f$ over $\mathbb C$ satisfying an equation of the
form~\eqref{caract theta quot bis}.
Using the functional equation~\eqref{func equa theta k}, we see that the~$k$-periodic meromorphic function
$g(z)=\frac{f(z) \theta_{k}(z-\xi)}{\theta_{k}(z)^{n} \theta_{k}(z)}$, where $\xi \in \mathbb C$ is such that $e^{-2i\pi
\xi/k}=(-1)^{n}c$, satisf\/ies $g(z+k\tau)=g(z)$.
So~$g$ belongs to $\operatorname{M}_{k\Lambda}^{\times} \subset \quotprodtheta{k}$, whence the result.
\end{proof}

\begin{Lemma}
If $f \in \prodtheta{k}$ is such that $\degr{k} (f)=0$ then~$f$ is constant.
\end{Lemma}

\begin{proof}
Consider $f \in \prodtheta{k}$.
There exists $c \in \mathbb C^{\times}$ and $(n_{\xi})_{\xi \in \mathbb C}\in \mathbb N^{(\mathbb C)}$ with f\/inite
support such that
\begin{gather*}
f(z)=c  \prod\limits_{\xi \in \mathbb C} \theta_{k}(z-\xi)^{n_{\xi}}.
\end{gather*}
Then, $\degr{k} (f)=\sum\limits_{\xi \in \mathbb C} n_{\xi}$ is equal to $0$ by hypothesis.
Thus, for all $\xi \in \mathbb C$, $n_{\xi}=0$ and hence $f=c$ is constant.
\end{proof}

\subsection[Weierstrass $\wp$-function]{Weierstrass $\boldsymbol{\wp}$-function}

For details on what follows, we refer to~\cite[Chapter~VI]{Silv}.
Recall that
\begin{gather*}
\wp(z):=\frac{1}{z^{2}}+ \sum\limits_{\lambda \in \Lambda \setminus \{0\}} \frac{1}{(z+\lambda)^{2}}
-\frac{1}{\lambda^{2}}\in \operatorname{M}_{\Lambda}
\end{gather*}
denotes the Weierstrass elliptic function associated to the lattice~$\Lambda$.
For any integer $k \geq 1$, we denote by $\wpf{k} \in \operatorname{M}_{k \Lambda}$ the Weierstrass function def\/ined by
\begin{gather*}
\wpf{k}(z):=\wp(z/k)\in \operatorname{M}_{k\Lambda}.
\end{gather*}
This $k\Lambda$-periodic meromorphic function is an even function, its poles are of order two and its set of poles is
$k\Lambda$.
Therefore, its derivative $\wpf{k}'$ is an odd function, its poles are of order three and its set of poles is
$k\Lambda$.

Any $k\Lambda$-periodic elliptic function is a~rational function in $\wpf{k}$ and $\wpf{k}'$, that is
\begin{gather*}
\operatorname{M}_{k\Lambda}=\mathbb C(\wpf{k},\wpf{k}').
\end{gather*}

\begin{Lemma}
\label{lem2}
Assume that $f \in \operatorname{M}_{k\Lambda}$, seen has a~meromorphic function over $\mathbb C/k\Lambda$, has at
most~$N$ poles counted with multiplicities $($or, equivalently, that $f=p/q$ with $p,q \in \prodtheta{k}$ such that
${\degr{k} p, \degr{k} q \leq N})$.
Then, there exist $A=P/Q$ and $B=R/S$ with~$P,Q \in \mathbb C[X]$ of degree at most $2N$ and $R,S \in \mathbb C[X]$ of
degree at most~$2N+3$ such that
\begin{gather*}
f=A(\wpf{k})+\wpf{k}' B(\wpf{k}).
\end{gather*}
\end{Lemma}

\begin{proof}
Using the fact that $f(z)$ belongs to $\operatorname{M}_{k\Lambda}$ if and only if $f(kz)$ belongs to
$\operatorname{M}_{\Lambda}$, it is easily seen that it is suf\/f\/icient to prove the lemma for $k=1$.
In what follows, we see the~$\Lambda$-periodic elliptic functions as meromorphic functions on $\mathbb C/\Lambda$.
Let~$A,B \in \mathbb C(X)$ be such that $f=A(\wp)+\wp' B(\wp)$.
It follows from the formula
\begin{gather*}
A(\wp(z))=\frac{f(z)+f(-z)}{2}
\end{gather*}
that $A(\wp)$ has at most $2N$ poles counted with multiplicities in $\mathbb C/\Lambda$.
But, if $A=P/Q$ with $\gcd(P,Q)=1$ then $A(\wp)$ has at least $\deg Q$ poles counted with multiplicities in $\mathbb
C/\Lambda$ (namely, the zeroes of $Q(\wp)$).
So $\deg Q \leq 2N$.

Using the fact that elliptic functions have the same numbers of zeroes and poles, the same argument applied to
$1/A(\wp)$ shows that~${\deg P \leq 2N}$.

Using the formula
\begin{gather*}
B(\wp(z))=\frac{f(z)-f(-z)}{2 \wp'(z)},
\end{gather*}
similar arguments show that $\deg R \leq 2N+3$ and $\deg S \leq 2N+3$.
\end{proof}

\section{Irreducibility of the dif\/ference Galois group}
\label{sec4}
We let
\begin{gather}
\label{equation d ordre 2}
\phi^{2}(y) + a \phi(y)+by=0 \qquad \text{with} \quad a \in \operatorname{M}_{\Lambda} \quad \text{and} \quad b\in
\operatorname{M}_{\Lambda}^{\times}
\end{gather}
be a~dif\/ference equation of order $2$ with coef\/f\/icients in $\operatorname{M}_{\Lambda}$ and we denote~by
\begin{gather*}
\phi Y=AY \qquad \hbox{with} \quad A=
\begin{pmatrix}
0&1
\\
-b&-a
\end{pmatrix}
\in \operatorname{GL}_{2}(\operatorname{M}_{\Lambda})
\end{gather*}
the associated dif\/ference system.
For the notations $\operatorname{M}_{\Lambda}$,~$\phi$,~$K$, etc, we refer to Sections~\ref{sec2} and~\ref{sec3}.

We let $G \subset \operatorname{GL}_{2}(\mathbb C)$ be the dif\/ference Galois group over $(K,\phi)$ of
equation~\eqref{equation d ordre 2}.
According to Corollary~\ref{coro1},~$G$ is an algebraic subgroup of $\operatorname{GL}_{2}(\mathbb C)$ such that the
quotient $G/G^{\circ}$ of~$G$ by its identity component $G^{\circ}$ is cyclic.
A~direct inspection of the classif\/ication, up to conjugation, of the algebraic subgroups of
$\operatorname{GL}_{2}(\mathbb C)$ given in~\cite[Theorem~4]{NvdPT} shows that~$G$ satisf\/ies one of the following
properties:
\begin{itemize}
\itemsep=0pt
\item The group~$G$ is reducible (i.e., conjugate to some subgroup of the group of upper-triangular matrices in
$\operatorname{GL}_{2}(\mathbb C)$).
If~$G$ is reducible, we distinguish the following sub-cases:
\begin{itemize}\itemsep=0pt
\item the group~$G$ is completely reducible (i.e., is conjugate to some subgroup of the group of diagonal matrices
in $\operatorname{GL}_{2}(\mathbb C)$);
\item the group~$G$ is not completely reducible.
\end{itemize}
\item The group~$G$ is irreducible (i.e., not reducible) and imprimitive (see Section~\ref{sec5} for the def\/inition).
\item The group~$G$ is irreducible and is not imprimitive, and, in this case, there exists an algebraic subgroup~$\mu$
of $\mathbb C^{\times}$ such that $G= \mu \operatorname{SL}_{2}(\mathbb C)$.
Therefore, $G= \{M \in \operatorname{GL}_{2}(\mathbb C) \, \vert \, \det(M) \in H\}$ where $H=\det(G) \subset \mathbb
C^{\times}$.
In order to determine~$H$, one can use the fact that $H=\det(G)$ is the dif\/ference Galois group of $\phi y = (\det A) y
= by$ (this follows for instance from Tannakian duality~\cite[Section~1.4]{VdPS97}).
\end{itemize}

Our f\/irst task, undertaken in the present section, is to study the reducibility of~$G$.
The imprimitivity of~$G$ will be considered in~Section~\ref{sec5}.

\subsection{Riccati equation and irreducibility}

The non linear dif\/ference equation
\begin{gather}
\label{eq4}
(\phi(u)+a) u=-b
\end{gather}
is called the Riccati equation associated to equation~\eqref{equation d ordre 2}.
A~straightforward calculation shows that~$u$ is a~solution of this equation if and only if $\phi-u$ is a~right factor of
$\phi^{2}+a\phi+b$, whence its link with irreducibility.

In what follows, we denote by $I_{2}$ the identity matrix of $\operatorname{GL}_2(\mathbb C)$.

\begin{Lemma}
\label{lem1}
The following statements hold:
\begin{enumerate}\itemsep=0pt
\item[$1.$] If~\eqref{eq4} has one and only one solution in~$K$ then~$G$ is reducible but not completely reducible.
\item[$2.$] If~\eqref{eq4} has exactly two solutions in~$K$ then~$G$ is completely reducible but not an algebraic subgroup of
$\mathbb C^{\times} I_{2}$.
\item[$3.$] If~\eqref{eq4} has at least three solutions in~$K$ then it has infinitely many solutions in~$K$ and~$G$ is an
algebraic subgroup of $\mathbb C^{\times} I_{2}$.
\item[$4.$] If none of the previous cases occur then~$G$ is irreducible.
\end{enumerate}
\end{Lemma}

\begin{proof}
The proof of this lemma is identical to that of~\cite[Theorem~4.2]{Hen98}, to whom we refer for more details.

(1) We assume that~\eqref{eq4} has one and only one solution $u\in K$.
A~straightforward calculation shows that
\begin{gather*}
\phi(T)AT^{-1}=
\begin{pmatrix}
u&*
\\
0&b/u
\end{pmatrix}
\qquad \text{for} \quad T:=
\begin{pmatrix}
1-u&1
\\
-u&1
\end{pmatrix}
\in\operatorname{GL}_{2}(K).
\end{gather*}
We deduce from this and from Corollary~\ref{coro1} that~$G$ is reducible.

Moreover, if~$G$ was completely reducible then, in virtue of Corollary~\ref{coro1}, $\phi(T)AT^{-1}$ would be diagonal
for some ${T:= (t_{i,j})_{1 \leq i,j\leq 2}\in \operatorname{GL}_{2}(K)}$.
Equating the entries of the antidiagonal of $\phi(T)AT^{-1}$ with~$0$, we f\/ind that
$-\frac{t_{21}}{t_{22}},-\frac{t_{11}}{t_{12}}\in K $ are solutions of the Riccati equation.
Since $\det(T) \neq 0$, these solutions are distinct, whence a~contradiction.

(2) Assume that~\eqref{eq4} has exactly two solutions $u_{1},u_{2}\in K$.
We have
\begin{gather*}
\phi(T)AT^{-1}=
\begin{pmatrix}
u_{1}&0
\\
0&u_{2}
\end{pmatrix}
\qquad \text{for}\quad T:=\frac{1}{u_{1}-u_{2}}
\begin{pmatrix}
-u_{2}&1
\\
-u_{1}&1
\end{pmatrix}
\in \operatorname{GL}_{2}(K).
\end{gather*}
We deduce from this and from Corollary~\ref{coro1} that~$G$ is completely reducible.

Moreover, if~$G$ was an algebraic subgroup of $\mathbb C^{\times} I_{2}$ then, according to Corollary~\ref{coro1}, there
would exist $u\in K$ and $T= (t_{i,j})_{1 \leq i,j\leq 2} \in \operatorname{GL}_{2}(K)$ such that
\begin{gather*}
\phi(T)AT^{-1}=u I_{2}.
\end{gather*}
This equality implies that $t_{21}$ and $t_{22}$ are non zero and that, for all $c,d\in \mathbb C$ with
$ct_{2,2}+dt_{1,2} \neq 0$,
\begin{gather*}
-\frac{ct_{21}+dt_{11}}{ct_{22}+dt_{12}}\in K
\end{gather*}
is a~solution of~\eqref{eq4}.
It is easily seen that we get in this way inf\/initely many solutions of the Riccati equation, this is a~contradiction.

(3) Assume that~\eqref{eq4} has at least three solutions $u_{1},u_{2},u_{3}\in K$.
The proof of assertion (2) of the present lemma shows that $\phi Y=AY$ is isomorphic over~$K$ to $\phi Y=
\begin{pmatrix}
u_{i}&0
\\
0&u_{j}
\end{pmatrix}
Y$ for all $1\leq i<j \leq 3$.
Therefore, there exists $T \in \operatorname{GL}_{2}(K)$ such that
\begin{gather*}
\phi(T)
\begin{pmatrix}
u_{1}&0
\\
0&u_{2}
\end{pmatrix}
=
\begin{pmatrix}
u_{1}&0
\\
0&u_{3}
\end{pmatrix}
T.
\end{gather*}
Equating the second columns in this equality, we see that there exists ${f\in K^{\times}}$ such that either
$u_{1}=\frac{\phi f}{f}u_{2}$ or $u_{3}=\frac{\phi f}{f}u_{2}$; up to renumbering, one can assume that the former case
holds true.
It follows that $\phi Y=AY$ is isomorphic over~$K$ to
\begin{gather*}
\phi Y=(u_{1} I_{2}) Y
\end{gather*}
and, according to Corollary~\ref{coro1},~$G$ is an algebraic subgroup of $\mathbb C^{\times} I_{2}$.
We have shown during the proof of statement (2) that this implies that the Riccati equation~\eqref{eq4} has inf\/initely
many solutions in~$K$.

(4) Assume that~$G$ is reducible.
According to Corollary~\ref{coro1}, there exists $T=(t_{i,j})_{1 \leq i,j\leq 2}\in \operatorname{GL}_{2}(K)$ such that
$\phi(T)AT^{-1}$ is upper triangular.
Then $t_{22}\neq 0$ and $-\frac{t_{21}}{t_{22}}\in K$ is a~solution of the Riccati equation~\eqref{eq4}.
This proves claim (4).
\end{proof}

In the proof of the previous lemma, we have shown the following result, which we state independently for ease of
reference.

\begin{Lemma}
\label{lem ric irred bis}
The following properties are equivalent:
\begin{itemize}\itemsep=0pt
\item The Riccati equation~\eqref{eq4} has at least three solutions in~$K$.
\item The Riccati equation~\eqref{eq4} has infinitely many solutions in~$K$.
\item The difference Galois group~$G$ is a~subgroup of $\mathbb C^{\times}I_{2}$.
\item There exist $u \in K^{\times}$ and $T \in \operatorname{GL}_2(K)$ such that $\phi(T)AT^{-1} = u I_{2}$.
\end{itemize}
\end{Lemma}

We shall now state and prove one more lemma.

\begin{Lemma}
\label{lem descente}
Let $\Lambda'' \subset \Lambda'$ be sublattices of~$\Lambda$ such that the quotient $\Lambda'/\Lambda''$ is cyclic.
Assume that there exist $u \in \operatorname{M}_{\Lambda''}^{\times}$ and $T \in
\operatorname{GL}_{2}(\operatorname{M}_{\Lambda''})$ such that
\begin{gather}
\label{et une equation une pour lem}
\phi(T)AT^{-1}=u I_{2}.
\end{gather}
Then, the Riccati equation~\eqref{eq4} has at least two distinct solutions in $\operatorname{M}_{\Lambda'}$.
\end{Lemma}

\begin{proof}
The Galois extension $ \operatorname{M}_{\Lambda''} \vert \operatorname{M}_{\Lambda'} $ is cyclic of order
$k:=[\operatorname{M}_{\Lambda''}: \operatorname{M}_{\Lambda'}]$.
Its Galois group $\operatorname{Gal}(\operatorname{M}_{\Lambda''}|\operatorname{M}_{\Lambda'})$ is generated by the
f\/ield automorphism $\sigma_{1}$ given by $\sigma_{1}(f(z))=f(z+\lambda')$, where $\lambda' \in \Lambda'$ is
a~representative of a~generator of $\Lambda'/\Lambda''$.
Note that the action of $\operatorname{Gal}(\operatorname{M}_{\Lambda''}|\operatorname{M}_{\Lambda'})$ on~$\operatorname{M}_{\Lambda''}$ commutes with the action of~$\phi$.
Applying~$\sigma_{1}$ to equation~\eqref{et une equation une pour lem}, we get
\begin{gather*}
\phi(\sigma_{1} (T))A\sigma_{1}(T)^{-1}=\sigma_{1}(u) I_{2},
\end{gather*}
so
\begin{gather*}
\phi(S)u=\sigma_{1}(u) S, \qquad \text{with} \quad S:=\sigma_{1}(T) T^{-1} \in \operatorname{GL}_{2}(\operatorname{M}_{\Lambda''}).
\end{gather*}
It follows that there exists $g_{\sigma_{1}}\in \operatorname{M}_{\Lambda''}^{\times}$ (namely, one of the non zero
entries of~$S$) such that
\begin{gather*}
\sigma_{1}(u)=\frac{\phi (g_{\sigma_{1}})}{g_{\sigma_{1}}}u.
\end{gather*}
Consider the norm
\begin{gather*}
\operatorname{N}:=\operatorname{N}_{\operatorname{M}_{\Lambda''}|\operatorname{M}_{\Lambda'}}(g_{\sigma_1})=
\prod\limits_{\sigma \in \operatorname{Gal}(\operatorname{M}_{\Lambda''}|\operatorname{M}_{\Lambda'})} \sigma
(g_{\sigma_1}) \in \operatorname{M}_{\Lambda'}^{\times}.
\end{gather*}
We have
\begin{gather*}
\phi (\operatorname{N}) = \prod\limits_{\sigma \in
\operatorname{Gal}(\operatorname{M}_{\Lambda''}|\operatorname{M}_{\Lambda'})} \sigma \left(\frac{\sigma_{1}(u)
g_{\sigma_{1}}}{u} \right) = \prod\limits_{\sigma \in
\operatorname{Gal}(\operatorname{M}_{\Lambda''}|\operatorname{M}_{\Lambda'})} \sigma(g_{\sigma_{1}}) = \operatorname{N},
\end{gather*}
so $\operatorname{N}=c \in (K^{\phi})^{\times} = \mathbb C^{\times}$.
Up to replacing $g_{\sigma_1}$ by $g_{\sigma_1}c^{-1/k}$, we may assume that
\begin{gather*}
\operatorname{N}_{\operatorname{M}_{\Lambda''}|\operatorname{M}_{\Lambda'}}(g_{\sigma_1})=1.
\end{gather*}
Hilbert's~90 theorem~\cite[Section~X.1]{Ser} ensures that there exists $m \in \operatorname{M}_{\Lambda''}^{\times}$ such
that
\begin{gather*}
{g_{\sigma_1}=\frac{m}{\sigma_1(m)}}.
\end{gather*}
For any $\sigma=\sigma_{1}^{j} \in \operatorname{Gal}(\operatorname{M}_{\Lambda''}|\operatorname{M}_{\Lambda'})$, we set
\begin{gather*}
g_\sigma:= g_{\sigma_1} \sigma_1 (g_{\sigma_1}) \cdots \sigma_1^{j-1} (g_{\sigma_1})=m/\sigma (m)\in
\operatorname{M}_{\Lambda''}^{\times};
\end{gather*}
we have
\begin{gather*}
\sigma(u)=\frac{\phi (g_{\sigma})}{g_{\sigma}}u.
\end{gather*}
It follows that
\begin{gather*}
\widetilde u:= \frac{\phi (m)}{m} u
\end{gather*}
is invariant under the action of $\operatorname{Gal}(\operatorname{M}_{\Lambda''}|\operatorname{M}_{\Lambda'})$ and
hence belongs to $\operatorname{M}_{\Lambda'}^{\times}$.
We have
\begin{gather*}
\phi\left(T'\right)A\left(T'\right)^{-1}= \widetilde{u} I_{2}, \qquad \text{with} \quad T':=mT \in
\operatorname{GL}_{2}(\operatorname{M}_{\Lambda''}).
\end{gather*}
Applying $\sigma \in \operatorname{Gal}(\operatorname{M}_{\Lambda''}|\operatorname{M}_{\Lambda'})$ to this equality, we
get
\begin{gather*}
\phi\left(\sigma(T')\right)A\left(\sigma(T')\right)^{-1}= \widetilde{u} I_{2}.
\end{gather*}
It follows that the matrix
\begin{gather*}
C_{\sigma}:=T' \sigma (T' )^{-1} \in \operatorname{GL}_{2}(\operatorname{M}_{\Lambda''})
\end{gather*}
satisf\/ies $\phi (C_{\sigma}) = C_{\sigma}$ and, hence, that its entries belong to~$K^{\phi} = \mathbb C$.
Moreover, $\sigma \mapsto C_{\sigma}$ is a~$1$-cocyle for the natural action of
$\operatorname{Gal}(\operatorname{M}_{\Lambda''}|\operatorname{M}_{\Lambda'})$ on $\operatorname{GL}_{2}(\mathbb C)$
but this action is trivial so $\sigma \mapsto C_{\sigma}$ is a~group morphism
from~$\operatorname{Gal}(\operatorname{M}_{\Lambda''}|\operatorname{M}_{\Lambda'})$ to $\operatorname{GL}_{2}(\mathbb
C)$.
Since $\operatorname{Gal}(\operatorname{M}_{\Lambda''}|\operatorname{M}_{\Lambda'})$ is cyclic, this implies that there
exists $P\in \operatorname{GL}_{2}(\mathbb C)$ such that, for all $\sigma\in
\operatorname{Gal}(\operatorname{M}_{\Lambda''}|\operatorname{M}_{\Lambda'})$, the matrix
${D_{\sigma}:=P^{-1}C_{\sigma}^{-1}P} \in \operatorname{GL}_{2}(\mathbb C)$ is diagonal.
We have, for all $\sigma \in\operatorname{Gal}(\operatorname{M}_{\Lambda''}|\operatorname{M}_{\Lambda'})$,
\begin{gather*}
\sigma(T'')=D_{\sigma}T'',\qquad  \text{where} \quad T''=(t''_{i,j})_{1 \leq i,j \leq 2}:=P^{-1}T' \in
\operatorname{GL}_{2}(\operatorname{M}_{k \Lambda}).
\end{gather*}
It follows that $u_{1}:=\frac{-t''_{11}}{t''_{12}}$ and $v_{1}:=\frac{-t''_{21}}{t''_{22}}$ are invariant by the action
of $\operatorname{Gal}(\operatorname{M}_{\Lambda''}|\operatorname{M}_{\Lambda'})$ and hence belong to
$\operatorname{M}_{\Lambda'}$.
But $u_{1}$ and $v_{1}$ are solutions of the Riccati equation~\eqref{eq4} (this was already used in the proof of
assertion (2) of Lemma~\ref{lem1}).
So $u_{1}$ and $v_{1}$ are solutions in $\operatorname{M}_{\Lambda'}$ of the Riccati equation~\eqref{eq4}.
\end{proof}

We now come to the main result of this subsection.

\begin{Theorem}
\label{theo2}
The following statements hold:
\begin{enumerate}
\itemsep=0pt
\item[$1.$] The Galois group~$G$ is reducible if and only if the Riccati equation~\eqref{eq4} has at least one solution
in~$\operatorname{M}_{2\Lambda}$.
\item[$2.$] The Galois group~$G$ is completely reducible if and only if the Riccati equation~\eqref{eq4} has at least two
solutions in~$\operatorname{M}_{2\Lambda}$.
\end{enumerate}
\end{Theorem}

\begin{proof}
In virtue of Lemma~\ref{lem1}, it is suf\/f\/icient to prove that:
\begin{enumerate}\itemsep=0pt
\item[(a)] If the Riccati equation~\eqref{eq4} has a~unique solution in~$K$, then it belongs to
$\operatorname{M}_{\Lambda}$.

\item[(b)] If the Riccati equation~\eqref{eq4} has exactly two solutions in~$K$, then they belong to
$\operatorname{M}_{2\Lambda}$.

\item[(c)] If the Riccati equation~\eqref{eq4} has at least three solutions in~$K$, then the Riccati equation~\eqref{eq4}
has at least two solutions in $\operatorname{M}_{2\Lambda}$.
\end{enumerate}

 (a) Assume that the Riccati equation~\eqref{eq4} has a~unique solution~$u$ in~$K$.
Since $u(z)$, $u(z+1)$ and $u(z+\tau)$ are solutions of~\eqref{eq4}, we get
\begin{gather*}
u(z)=u(z+1)=u(z+\tau)
\end{gather*}
and hence $u\in\operatorname{M}_{\Lambda}$.

(b) Assume that the Riccati equation~\eqref{eq4} has exactly two solutions in~$K$ and let $u\in K$ be one of these
solutions.
Since $u(z)$, $u(z+1)$ and $u(z+2)$ are solutions of~\eqref{eq4}, we get $ u(z+2)=u(z).
$ Similarly, we have $u(z+2\tau)=u(z)$.
So $u\in\operatorname{M}_{2\Lambda}$.

(c) What follows is inspired by~\cite[Theorem~4.2]{Hen98}, but is a~little bit subtler.
Assume that the Riccati equation~\eqref{eq4} has at least three solutions in~$K$.
According to Lemma~\ref{lem ric irred bis}, there exist $u\in K$ and ${T=(t_{i,j})_{1 \leq i,j \leq 2} \in
\operatorname{GL}_{2}(K)}$ such that
\begin{gather}
\label{et une equation une}
\phi(T)AT^{-1}=u I_{2}.
\end{gather}
Let $k \in \mathbb N^{*}$ be such that the entries of~$T$ and~$u$ belong to $\operatorname{M}_{k \Lambda}$.
Consider the following f\/ield extensions:
\begin{gather*}
\operatorname{M}_{\Lambda} \subset L \subset \operatorname{M}_{k\Lambda}, \qquad \text{with}\quad L:=\operatorname{M}_{\mathbb
Z+k\tau\mathbb Z}.
\end{gather*}
Applying Lemma~\ref{lem descente} to the extension $\operatorname{M}_{k\Lambda}\vert L$ and to the equation~\eqref{et
une equation une}, we get that the Riccati equation~\eqref{eq4} has two distinct solution $u_{1}$ and $v_{1}$ in~$L$.
If both of them belong to $\operatorname{M}_{2\Lambda}$ then the proof is completed.
Otherwise, up to renumbering, we can assume that $u_{1} \not \in \operatorname{M}_{2\Lambda}$, i.e., that $u_{1}$ is
not $2\tau$-periodic.
Then
\begin{gather*}
u_{1}(z),u_{2}(z):=u_{1}(z+\tau) \qquad \text{and} \qquad u_{3}(z):=u_{1}(z+2\tau)
\end{gather*}
are distinct solutions in~$L$ of the Riccati equation.
For all integers $i,j \in \{1,2,3\}$ with $i<j$ we set $T_{i,j}:=\frac{1}{u_{i}-u_{j}}
\begin{pmatrix}
-u_{j}&1
\\
-u_{i}&1
\end{pmatrix}
\in \operatorname{GL}_{2}(L)$ and we have
\begin{gather*}
\phi (T_{i,j} )A (T_{i,j} )^{-1}=
\begin{pmatrix}
u_{i}&0
\\
0&u_{j}
\end{pmatrix}
\end{gather*}
(this was already used in the proof of assertion~(2) of Lemma~\ref{lem1}).
Therefore,
\begin{gather*}
\phi\big(T_{1,3} (T_{1,2} )^{-1}\big)
\begin{pmatrix}
u_{1}&0
\\
0&u_{2}
\end{pmatrix}
=
\begin{pmatrix}
u_{1}&0
\\
0&u_{3}
\end{pmatrix}
T_{1,3} (T_{1,2} )^{-1}.
\end{gather*}
Equating the second columns in this equality, we see that there exists $f \in L^{\times}$ such that either
$u_{1}=\frac{\phi f}{f} u_{2}$ or $u_{3}=\frac{\phi f}{f} u_{2}$; up to renumbering, we may assume that the former
equality holds true.
Then, we have
\begin{gather}
\label{et une deux equa}
\phi\big(\widetilde{T}\big)A\widetilde{T}^{-1}= u_{1} I_{2}
\end{gather}
with
\begin{gather*}
u_{1} \in L^{\times} \qquad \text{and} \qquad \widetilde{T}:=
\begin{pmatrix}
1&0
\\
0&f
\end{pmatrix}
T_{1,2}\in \operatorname{GL}_{2}(L).
\end{gather*}
Applying Lemma~\ref{lem descente} to the extension $L \vert \operatorname{M}_{\Lambda}$ and to the equation~\eqref{et
une deux equa}, we see that the Riccati equation~\eqref{eq4} has $2$ distinct solutions in $\operatorname{M}_{\Lambda}$.
This concludes the proof.
\end{proof}

\subsection{On the solutions of the Riccati equation}

We refer to Section~\ref{sec31} for the notations ($\divi{k}$, $\degr{k}$, $\weight{k}$, etc.) used in this subsection.
Let $k\geq 1$ be an integer.
Consider $p_{1} \in \prodtheta{k} \cup \{0\}$ and $p_{2},p_{3} \in \prodtheta{k}$ such that
\begin{gather*}
a= \frac{p_{1}}{p_{3}} \qquad  \text{and} \qquad b = \frac{p_{2}}{p_{3}}.
\end{gather*}
We let $u \in \operatorname{M}_{k \Lambda}$ be a~potential solution of the Riccati equation~\eqref{eq4}.

\begin{Proposition}
\label{propo2}
We have
\begin{gather*}
u = \frac{\phi(r)}{r}\frac{p}{q}
\end{gather*}
for some $p,q, r\in \prodtheta{k}$ such that
\begin{itemize}
\itemsep=0pt
\item[$(i)$] $\divi{k} (p) \leq \divi{k} (p_{2})$,
\item[$(ii)$] $\divi{k} (q) \leq \divi{k} (\phi^{-1}(p_{3}))$,
\item[$(iii)$] $\degr{k} (p)=\degr{k} (q)$,
\item[$(iv)$] $\weight{k} (p/q)=\degr{k} (r) h \mod k\Lambda$.
\end{itemize}
\end{Proposition}

\begin{proof}
In what follows, the greatest common divisors $(\gcd)$ has to be understood in the ring~$\mathcal{O}(\mathbb C)$ of
entire functions\footnote{According to~\cite{Helmer}, any f\/initely generated ideal of $\mathcal{O}(\mathbb C)$ is
principal, whence the existence of the greatest common divisor of any couple of elements of $\mathcal{O}(\mathbb C)$.
Such a~greatest common divisor is unique up to multiplication by an unit of $\mathcal{O}(\mathbb C)$.}.
Let $p_{4}, p_{5} \in \prodtheta{k}$, with $\gcd(p_{4},p_{5})=1$, be such that $u=p_{4}/p_{5}$.
Let $r \in \prodtheta{k}$ be a~greatest common divisor of $\phi^{-1}(p_{4})$ and $p_{5}$ and consider
\begin{gather*}
p:=\frac{p_{4}}{\phi(r)} \in \prodtheta{k} \qquad \text{and} \qquad q:=\frac{p_{5}}{r} \in \prodtheta{k}.
\end{gather*}
By construction, we have
\begin{gather*}
u= \frac{\phi r}{r} \frac{p}{q}
\end{gather*}
with $\gcd(p,\phi (q))=\gcd(\phi (r) p,r q)=1$.
Then, the Riccati equation~\eqref{eq4} becomes
\begin{gather*}
p_{3} \frac{\phi r}{r} \frac{p}{q} \phi \left(\frac{\phi r}{r} \frac{p}{q}\right)+p_{1} \frac{\phi r}{r}
\frac{p}{q}=-p_{2},
\end{gather*}
i.e.,
\begin{gather*}
p_{3} \phi^{2}(r) p \phi (p)+p_{1} \phi(r) p \phi(q)=-p_{2} r q\phi(q).
\end{gather*}
It is now easily seen that~$p$ divides $p_{2}$ and that~$q$ divides $\phi^{-1}(p_{3})$ in $\mathcal{O}(\mathbb C)$.
In terms of divisors, this is exactly (i) and (ii).

According to Lemma~\ref{form precise}, we have
\begin{gather*}
\frac{p}{q}(z+k\tau)= (-1)^{\degr{k} (p/q)} e^{2i\pi \omega/k} e^{-2 i\pi \degr{k} (p/q) z/k} \frac{p}{q}(z)
\end{gather*}
for some representative~$\omega$ of $\weight{k}(p/q)$, and
\begin{gather*}
\frac{\phi(r)}{r}(z+k\tau)=e^{- 2 i\pi \degr{k} (r) h/k} \frac{\phi(r)}{r}(z).
\end{gather*}
Therefore
\begin{gather*}
u(z+k\tau)= (-1)^{\degr{k} (p/q)} e^{2i\pi \omega/k} e^{-2 i\pi \degr{k} (p/q) z/k} e^{- 2 i\pi \degr{k} (r) h/k} u(z).
\end{gather*}
But $u \in \operatorname{M}_{k\Lambda}$, so $u(z+k\tau)=u(z)$ and, hence,
\begin{gather*}
(-1)^{\degr{k} (p/q)} e^{2i\pi \omega/k} e^{-2 i\pi \degr{k} (p/q) z/k} e^{- 2 i\pi \degr{k} (r) h/k}=1.
\end{gather*}
Hence $ \degr{k} (p/q)=0$ and $\omega = \degr{k} (r) h \mod k\Lambda$.
This proves~(iii) and~(iv).
\end{proof}

We will see in Section~\ref{sec61} that Proposition~\ref{propo2} is a~useful theoretic tool in order to determine the
dif\/ference Galois groups of families of equations, such as the discrete Lam\'e equations mentioned in the introduction.

We shall now conclude this section with a~few words about Proposition~\ref{propo2}.

\begin{Remark}
How to use Proposition~\ref{propo2} in order to decide whether~$G$ is irreducible? Theo\-rem~\ref{theo2} ensures that~$G$
is irreducible if and only if the Riccati equation~\eqref{eq4} has a~solution $u \in \operatorname{M}_{2 \Lambda}$; we
let~$p$,~$q$,~$r$ be as in Proposition~\ref{propo2}.
Assertions~(i) and~(ii) of Proposition~\ref{propo2}, show that there are \textit{finitely} many explicit possibilities
for the divisors $\divi{2} (p)$ and~$\divi{2}(q)$.
But $\degr{2} (r)$ is entirely determined by these divisors in virtue of~(iv) of Proposition~\ref{propo2}.
So, we can compute an integer $N \geq 0$ such that if the Riccati equation~\eqref{eq4} has a~solution $u \in
\operatorname{M}_{2 \Lambda}$, then
\begin{gather*}
u=p_{0}/q_{0}
\end{gather*}
with $p_{0},q_{0} \in \prodtheta{2}$ such that~${\degr{2} (p_{0}) \leq N}$ and ${\degr{2} (q_{0}) \leq N}$.
Lemma~\ref{lem2} ensures that
\begin{gather*}
u=A(\wpf{2})+\wpf{2}' B(\wpf{2})
\end{gather*}
for some $A=P/Q$ and $B=R/S$ with $P,Q \in \mathbb C[X]$ of degree at most $2N$ and $R,S \in \mathbb C[X]$ of degree at
most $2N+3$.

So, in order to determine whether or not the Riccati equation~\eqref{eq4} has at least one solution in
$\operatorname{M}_{2\Lambda}$, we are lead to the following question: do there exist $A=P/Q$ and $B=R/S$ with $P,Q \in
\mathbb C[X]$ of degree at most $2N$ and ${R,S \in \mathbb C[X]}$ of degree at most $2N+3$ such that $
u=A(\wpf{2})+\wpf{2}' B(\wpf{2}) $ is a~solution of the Riccati equation~\eqref{eq4}? Substituting
$u=A(\wpf{2})+\wpf{2}' B(\wpf{2})$ in the Riccati equation~\eqref{eq4} and using the addition formula:
\begin{gather*}
\wp_{2}(z)+\wp_{2}(h)+\wp_{2}(z+h)=\frac{1}{4} \left(\frac{\wp_{2}'(z)-\wp_{2}'(h)}{\wp_{2}(z)-\wp_{2}(h)}\right)^{2},
\end{gather*}
we are lead to decide whether multivariate polynomials, whose indeterminates are the coef\/f\/icients of~$P$,~$Q$,~$R$ and~$S$,
have a~common complex solution.
This can be decided by using Gr\"obner bases.

Note however that, in order to make this method an ef\/fective tool, we have to know the divisors of~$a$ and~$b$, and to
be able to deduce $\degr{k}(r)$ from assertion (iv) of Proposition~\ref{propo2}.
\end{Remark}

\section{Imprimitivity of the dif\/ference Galois group}
\label{sec5}

We want to determine whether~$G$ is imprimitive, that is whether~$G$ is conjugate to a~subgroup of
\begin{gather*}
\left\{
\begin{pmatrix}
\alpha&0
\\
0&\beta
\end{pmatrix}
\, \vert \, \alpha, \beta \in \mathbb C^{\times} \right\} \bigcup \left\{
\begin{pmatrix}
0&\gamma
\\
\delta &0
\end{pmatrix}
\, \vert \, \gamma, \delta \in \mathbb C^{\times} \right\}.
\end{gather*}

\begin{Theorem}
\label{riccati et imprimitivite}
Assume that~$G$ is irreducible and that $a\neq 0$.
Then,~$G$ is imprimitive if and only if there exists $u\in \operatorname{M}_{2\Lambda}$ such that
\begin{gather}
\label{eq6}
\left(\phi^{2}(u)+\left(\phi^{2}\left(\frac{b}{a}\right)-\phi(a)+\frac{\phi(b)}{a} \right)\right)
u=-\frac{\phi(b)b}{a^{2}}.
\end{gather}
\end{Theorem}

\begin{proof}
Arguing exactly as in~\cite[Theorem 4.6]{Hen98}, we get that~$G$ is imprimitive if and only if equation~\eqref{eq6} has
a~solution in~$K$.
But this is a~Riccati-type equation, with~$\phi$ replaced by $\phi^{2}$.
Therefore, the assertions (a), (b) and (c) given at the beginning of the proof of Theorem~\ref{theo2} allow us to
conclude.
\end{proof}

\begin{Remark}
If $a=0$ then~$G$ is imprimitive in virtue of Corollary~\ref{coro1}.
\end{Remark}

Note that Proposition~\ref{propo2} can be used in order to f\/ind restrictions on the solutions of the above Riccati-type
equation, but with~$\phi$ replaced by $\phi^{2}$.

\section{Applications}
\label{sec6}

We recall that $h \in \mathbb C$ is such that $h \mod \Lambda$ is not a~torsion point of $\mathbb C/\Lambda$, i.e.,
that the corresponding point $\underline{h}$ of $\mathcal{E}(\mathbb C)$ is not a~torsion point.

\subsection{A discrete version of Lam\'e equation}
\label{sec61}
Let us consider the dif\/ference equation
\begin{gather}
\label{lame deform}
\Delta_h^2 y = (A\wp(z)+B)y, \qquad  \text{where} \quad \Delta_h y(z)= \frac{y(z+h)-y(z)}{h}
\end{gather}
and $A,B \in \mathbb C$.
This is a~discrete version of the so-called Lam\'e dif\/ferential equation
\begin{gather*}
y''(z) = (A\wp(z)+B)y(z).
\end{gather*}

\begin{Theorem}
\label{theo lame alg}
Assume that $\mathcal{E}$ is defined over $\overline{\mathbb{Q}}$ $($i.e., $g_{2},g_{3} \in \overline{\mathbb{Q}})$
and that $h,A,B\in \overline{\mathbb{Q}}$ with $A \neq 0$.
Then, the difference Galois group over $(K,\phi)$ of equation~\eqref{lame deform} is $\operatorname{GL}_{2}(\mathbb C)$.
\end{Theorem}

A straightforward calculation shows that equation~\eqref{lame deform} can be rewritten as follows:
\begin{gather*}
\phi^2 y -2\phi y+\big({-}Ah^2\wp(z)-Bh^2+1\big)y=0.
\end{gather*}
We will deduce Theorem~\ref{theo lame alg} from the following theorem combined with a~transcendence result due to
Schneider.

\begin{Theorem}
\label{theo3}
Consider $a\in \mathbb C^{\times}$ and $b(z)=\alpha \wp(z) + \beta$ with $(\alpha,\beta) \in \mathbb C^{\times} \times
\mathbb C$.
Let $z_0 \in \mathbb C$ be such that $\wp(z_0)=-\beta/\alpha$.\footnote{Any non constant elliptic function $f(z)$ has at
least one zero (otherwise, $1/f(z)$ would be an entire elliptic function and hence would be constant).
In particular, $\wp(z)+\beta/\alpha$ has a~least one zero in $\mathbb C$.}
If $\mathbb Z h \cap (\ell z_{0} + \Lambda)=\{0\}$ for all $\ell \in \{-8,\dots,8\}$ $($this holds in particular if
$\mathbb Z h \cap (\mathbb Z z_{0} + \Lambda)=\{0\})$ then the difference Galois group over $(K,\phi)$ of $\phi^2 y
+a\phi y+by=0$ is $\operatorname{GL}_{2}(\mathbb C)$.
\end{Theorem}

\begin{proof}[Proof of Theorem~\ref{theo3}] For the notations, $\divi{k}$, $[\cdot]_k$, etc, we refer to Section~\ref{sec31}.
Note that
\begin{gather*}
\divi{1} (b)=[z_{0}]_{1}+[-z_{0}]_{1}-2[0]_{1}.
\end{gather*}
So, we can write $a=\frac{p_{1}}{p_{3}}$ and $b=\frac{p_{2}}{p_{3}}$ for some $p_{1},p_{2},p_{3}\in \prodtheta{1}$ with
\begin{gather*}
\divi{1} (p_{2})=[z_{0}]_{1}+[-z_{0}]_{1}
\end{gather*}
and
\begin{gather*}
\divi{1} (p_{3})=2[0]_{1}.
\end{gather*}

We claim that~$G$ is irreducible, i.e., in virtue of Theorem~\ref{theo2}, that the Riccati equation
\begin{gather}
\label{riccati irred lame}
(\phi(u)+a)u=-b
\end{gather}
does not have any solution in $\operatorname{M}_{2\Lambda}$.
Suppose to the contrary that it has a~solution $u \in \operatorname{M}_{2\Lambda}$.
Proposition~\ref{propo2} ensures that there exist $p,q, r\in \prodtheta{2}$ such that
\begin{gather*}
u= \frac{\phi(r)}{r}\frac{p}{q}
\end{gather*}
and
\begin{itemize}
\itemsep=0pt
\item[(i)] $\divi{2} (p) \leq \sum\limits_{\ell_{1},\ell_{2} \in
\{0,1\}}[\ell_{1}+\ell_{2}\tau-z_{0}]_{2}+[\ell_{1}+\ell_{2}\tau+z_{0}]_{2}$,
\item[(ii)] $\divi{2} (q) \leq \sum\limits_{\ell_{1},\ell_{2} \in \{0,1\}}2[\ell_{1}+\ell_{2}\tau+h]_{2}$,
\item[(iii)] $\degr{2} (p)=\degr{2} (q)$,
\item[(iv)] $\weight{2} (p/q)=\degr{2} (r) h \mod 2\Lambda$.
\end{itemize}
Properties (i) and (ii) above imply that
\begin{gather*}
\weight{2} \left(p/q\right) = \ell z_{0}-\degr{2} (q)h \mod \Lambda
\end{gather*}
for some $\ell \in \{-4,\dots,4\}$.
We infer from this and from (iv) that
\begin{gather*}
(\degr{2} (r) + \degr{2} (q))h = \ell z_{0} \mod \Lambda.
\end{gather*}
The assumption on $z_{0}$ ensures that $\degr{2} (r) = \degr{2} (q)=0$.
It follows from~(iii) that $ \degr{2} (p)=0$ and hence~$u$ is a~constant.
But it is easily seen that equation~\eqref{riccati irred lame} does not have any constant solution; this proves our
claim.

We claim that~$G$ is not imprimitive, i.e., in virtue of Theorem~\ref{riccati et imprimitivite}, that
\begin{gather}
\label{riccati imprim lame}
\left(\phi^{2}(u)+\frac{\phi^{2}(b)}{a}-a+\frac{\phi(b)}{a} \right)u=-\frac{\phi(b)b}{a^{2}}
\end{gather}
does not have any solution in $\operatorname{M}_{2\Lambda}$.
Suppose to the contrary that it has a~solution $u \in \operatorname{M}_{2\Lambda}$.
Equation~\eqref{riccati imprim lame} is of the form:
\begin{gather*}
u\left(\phi^{2}(u)+\frac{p_{1}}{p_{3}}\right)=\frac{p_{2}}{p_{3}},
\end{gather*}
for some $p_{1},p_{2},p_{3}\in \prodtheta{1}$ with
\begin{gather*}
\divi{1} (p_{2})=2[-2h]_{1}+[z_{0}]_{1}+[-z_{0}]_{1}+[z_{0}-h]_{1}+[-z_{0}-h]_{1}
\end{gather*}
and
\begin{gather*}
\divi{1} (p_{3})=2[-2h]_{1}+2[-h]_{1}+2[0]_{1}.
\end{gather*}
We apply Proposition~\ref{propo2} with~$\phi$ replaced by $\phi^2$ to obtain the existence of $p,q,r\in \prodtheta{2}$
such that
\begin{gather*}
u= \frac{\phi^{2}(r)}{r}\frac{p}{q},
\end{gather*}
where
\begin{itemize}\itemsep=0pt
\item[(v)]
\begin{gather*}
 \divi{2} (p)  \leq \sum\limits_{\ell_{1},\ell_{2} \in
\{0,1\}} 2[\ell_{1}+\ell_{2}\tau-2h]_{2}+[\ell_{1}+\ell_{2}\tau+z_{0}]_{2}+[\ell_{1}+\ell_{2}\tau-z_{0}]_{2}
\\
\hphantom{\divi{2} (p)  \leq}{} +[\ell_{1}+\ell_{2}\tau+z_{0}-h]_{2}+[\ell_{1}+\ell_{2}\tau-z_{0}-h]_{2},
\end{gather*}
\item[(vi)] $\divi{2} (q) \leq \sum\limits_{\ell_{1},\ell_{2} \in
\{0,1\}}2[\ell_{1}+\ell_{2}\tau]_{2}+2[\ell_{1}+\ell_{2}\tau+h]_{2}+2[\ell_{1}+\ell_{2}\tau+2h]_{2}$,
\item[(vii)] $\degr{2} (p)=\degr{2} (q)$,
\item[(viii)] $\weight{2} (p/q)=2\degr{2} (r) h \mod 2\Lambda$.
\end{itemize}
We claim that
\begin{itemize}
\itemsep=0pt
\item[(v$'$)] $\divi{2} (p) \leq \sum\limits_{\ell_{1},\ell_{2} \in
\{0,1\}}[\ell_{1}+\ell_{2}\tau+z_{0}]_{2}+[\ell_{1}+\ell_{2}\tau-z_{0}]_{2}$,
\item[(vi$'$)] $\divi{2} (q) \leq \sum\limits_{\ell_{1},\ell_{2} \in \{0,1\}}2[\ell_{1}+\ell_{2}\tau]_{2}$.
\end{itemize}
Indeed, otherwise, arguing as for the proof of the irreducibility of~$G$, we see that (v), (vi) and (viii) would lead to
a~relation of the form
\begin{gather*}
(2\degr{2} (r) + d)h = \ell z_{0} \mod \Lambda
\end{gather*}
for some integer $\ell \in \{-8,\dots,8\}$ and some integer $d>0$ and this would contradict our assumption on $z_{0}$.
Then, (viii) shows that
\begin{gather*}
2\degr{2} (r) h = \ell z_{0} \mod \Lambda
\end{gather*}
for some integer $\ell \in \{-4,\dots,4\}$ and hence $\degr{2} (r)=0$.
Therefore, $u=p/q$ with $p,q \in \prodtheta{2}$ satisfying (v$'$) and (vi$'$) above.
Now remark that
\begin{gather*}
\phi^{2}(u)+\frac{\phi^{2}(b)}{a}-a+\frac{\phi(b)}{a}
\end{gather*}
does not have poles in~$\Lambda$.
But any element of~$\Lambda$ is a~pole of order $2$ of the right hand side of equation~\eqref{riccati imprim lame}, so
any element of~$\Lambda$ is a~pole of order at least~$2$ of~$u$.
It follows that~(vi$'$) is an equality.
Then, using~(vii), we see that (v$'$) is also an equality.

So $\divi{2}(u)=\divi{2}(b)$ and hence $u=cb$ for some $c \in \mathbb C^{\times}$.
We now plug $u=cb$ into equation~\eqref{riccati imprim lame} and we get:
\begin{gather*}
c\left(\left(c+\frac{1}{a}\right)\phi^{2}(b)-a+\frac{\phi(b)}{a} \right)=-\frac{\phi(b)}{a^{2}}.
\end{gather*}
Since $-2h$ is a~pole of $\phi^{2}(b)$ but not of $\phi(b)$, we get $c=-1/a$ and the above equation simplif\/ies as
follows:
\begin{gather*}
\frac{-1}{a}\left(-a+\frac{\phi(b)}{a} \right)=-\frac{\phi(b)}{a^{2}}.
\end{gather*}
This gives $1=0$, whence a~contradiction.

Therefore,~$G$ is irreducible and not imprimitive.
So, as explained at the beginning of Section~\ref{sec4}, $G = \{M \in \operatorname{GL}_{2}(\mathbb C) \, \vert \, \det
(M) \in H\}$ where~$H \subset \mathbb C^{\times}$ is the Galois group of $\phi y = b y$, which is easily seen to be the
multiplicative group $(\mathbb C^{\times},\cdot)$.
This concludes the proof.
\end{proof}

\begin{proof}[Proof of Theorem~\ref{theo lame alg}]
In virtue of Theorem~\ref{theo3}, it is suf\/f\/icient to prove that $\mathbb Z h
\cap (\mathbb Z z_{0} + \Lambda)=\{0\}$.
Consider $m_{1},m_{2}\in \mathbb Z$ and $\lambda \in \Lambda$ such that $m_{1}h=m_{2}z_{0}+\lambda$.
We have $\wp(z_{0})=\frac{-Bh^{2}+1}{A h^{2}} \in \overline{\mathbb{Q}}$.
It follows that either $m_{2}z_{0}+\lambda \in \Lambda$ or $\wp(m_{2}z_{0}+\lambda) \in \overline{\mathbb{Q}}$.
(Indeed, suppose that $m_{2}z_{0}+\lambda \not \in \Lambda$.
Using equation~\eqref{Weierstrass equation}, we see that $\wp'(z_{0}) \in \overline{\mathbb{Q}}$.
Therefore, $\varphi(z_{0})$ belongs to $ \mathcal{E}(\overline{\mathbb{Q}})$, the map~$\varphi$ being def\/ined in the
introduction.
Using the fact that~$\varphi$ is a~group morphism and that $\mathcal{E}(\overline{\mathbb{Q}})$ is a~subgroup of
$\mathcal{E}(\mathbb C)$, we get $\varphi(mz_{0}) \in \mathcal{E}(\overline{\mathbb{Q}})$.
Therefore, $\wp(mz_{0}+\lambda)=\wp(mz_{0}) \in \overline{\mathbb{Q}}$.) In the former case, we get $m_{1}h\in \Lambda$
and hence $m_{1}=0$.
In the later case, it follows from the work of Schneider~\cite{schneider} (for a~reference in english, see Baker's
book~\cite[Theorem 6.2]{BakerTNT}; see also Bertrand and Masser's papers~\cite{BerMas,Mas75}) that $m_{2}z_{0}+\lambda$
and hence $m_{1}h$ are transcendental numbers, which is excluded.
\end{proof}

\subsection[A family of examples with Galois groups between $\operatorname{SL}
_{2}(\mathbb C)$ and~$\operatorname{GL}_{2}(\mathbb C)$]{A family of examples with Galois groups between $\boldsymbol{\operatorname{SL}
_{2}(\mathbb C)}$ and~$\boldsymbol{\operatorname{GL}_{2}(\mathbb C)}$}

\begin{Theorem}
\label{theo je ne sais plus combien}
Let us consider~$b\in \mathbb C^{\times}$, and $a(z):=\alpha \wp(z) + \beta$ with ${(\alpha,\beta) \in \mathbb
C^{\times}\times \mathbb C}$.
Let $z_0 \in \mathbb C$ be such that $\wp(z_0)=-\beta/\alpha$.
If $\mathbb Z h \cap (\ell z_{0} + \Lambda)=\{0\}$ for all $\ell \in \{-16,\dots,16\}$ $($this holds in particular if
$\mathbb Z h \cap (\mathbb Z z_{0} + \Lambda)=\{0\})$ then the difference Galois group over $(K,\phi)$ of $\phi^2 y
+a\phi y+by=0$ is $\mu_{2k} \operatorname{SL}_{2}(\mathbb C)$ if~$b$ is a~primitive~$k$th root of the unity and
$\operatorname{GL}_{2}(\mathbb C)$ otherwise, where $\mu_{2k}$ is the group of complex~$k$th roots of the unity.
\end{Theorem}

The proof will be given after the following corollary.

\begin{Corollary}
Assume that $\mathcal{E}$ is defined over $\overline{\mathbb{Q}}$ $($i.e., $g_{2},g_{3} \in \overline{\mathbb{Q}})$.
Consider $b\in \overline{\mathbb{Q}}^{\times}$ and $a(z):=\alpha \wp(z) + \beta$ with $\alpha,\beta \in
\overline{\mathbb{Q}}$ and $\alpha \neq 0$.
Then, the difference Galois group over $(K,\phi)$ of $\phi^2 y +a\phi y+by=0$ is $\mu_{2k} \operatorname{SL}_{2}(\mathbb
C)$ if~$b$ is a~primitive~$k$th root of the unity and $\operatorname{GL}_{2}(\mathbb C)$ otherwise.
\end{Corollary}

\begin{proof}
Similar to deduction of Theorem~\ref{theo lame alg} from Theorem~\ref{theo3}.
\end{proof}

\begin{proof}[Proof of Theorem~\ref{theo je ne sais plus combien}] Note that
\begin{gather*}
\divi{1} (a)=[z_{0}]_{1}+[-z_{0}]_{1}-2[0]_{1}.
\end{gather*}
So, we can write $a=\frac{p_{1}}{p_{3}}$ and $b=\frac{p_{2}}{p_{3}}$ for some $p_{1},p_{2},p_{3}\in \prodtheta{1}$ with
\begin{gather*}
\divi{1} (p_{2})=2[0]_{1}
\end{gather*}
and
\begin{gather*}
\divi{1} (p_{3})=2[0]_{1}.
\end{gather*}

We claim that~$G$ is irreducible, i.e., in virtue of Theorem~\ref{theo2}, that the Riccati equation
\begin{gather}
\label{eq2}
(\phi(u)+a)u=-b
\end{gather}
does not have any solution in $\operatorname{M}_{2\Lambda}$.
Suppose to the contrary that it has a~solution $u \in \operatorname{M}_{2\Lambda}$.
Proposition~\ref{propo2} ensures that there exist $p,q, r\in \prodtheta{2}$ such that
\begin{gather*}
u= \frac{\phi(r)}{r}\frac{p}{q}
\end{gather*}
and
\begin{itemize}
\itemsep=0pt
\item[(i)] $\divi{2} (p) \leq \sum\limits_{\ell_{1},\ell_{2} \in \{0,1\}}2[\ell_{1}+\ell_{2}\tau]_{2}$,
\item[(ii)] $\divi{2} (q) \leq \sum\limits_{\ell_{1},\ell_{2} \in \{0,1\}}2[\ell_{1}+\ell_{2}\tau+h]_{2}$,
\item[(iii)] $\degr{2} (p)=\degr{2} (q)$,
\item[(iv)] $\weight{2} (p/q)=\degr{2} (r) h \mod 2\Lambda$.
\end{itemize}
Properties (i) and (ii) above imply that
\begin{gather*}
\weight{2} \left(p/q\right) =-h\degr{2} (q) \mod \Lambda.
\end{gather*}
We infer from this and from (iv) that
\begin{gather*}
(\degr{2} (r) + \degr{2} (q))h =0 \mod \Lambda.
\end{gather*}
This yields $\degr{2} (r) = \degr{2} (q)=0$.
It follows from~(iii) that $ \degr{2} (p)=0$ and hence~$u$ is a~constant.
But it is easily seen that equation~\eqref{eq2} does not have any constant solution; this proves our claim.

We claim that~$G$ is not imprimitive, i.e., in virtue of Theorem~\ref{riccati et imprimitivite}, that (we recall
that~$b$ is constant)
\begin{gather}
\label{eq10}
\left(\phi^{2}(u)+\frac{b}{\phi^{2}(a)}-\phi(a)+\frac{b}{a} \right)u=-\frac{b^{2}}{a^{2}}
\end{gather}
does not have any solution in $\operatorname{M}_{2\Lambda}$.
Suppose to the contrary that it has a~solution $u \in \operatorname{M}_{2\Lambda}$.
Equation~\eqref{eq10} is of the form:
\begin{gather*}
u\left(\phi^{2}(u)+\frac{p_{1}}{p_{3}}\right)=\frac{p_{2}}{p_{3}},
\end{gather*}
for some $p_{1},p_{2},p_{3}\in \prodtheta{1}$ with
\begin{gather*}
\divi{1} (p_{2})=4[0]_{1}+2[-h]_{1}+[z_{0}-2h]_{1}+[-z_{0}-2h]_{1}
\end{gather*}
and
\begin{gather*}
\divi{1} (p_{3})=2[-h]_{1}+2[z_{0}]_{1}+2[-z_{0}]_{1}+[z_{0}-2h]_{1}+[-z_{0}-2h]_{1}.
\end{gather*}
Proposition~\ref{propo2} ensures that there exist $p,q,r\in \prodtheta{2}$ such that
\begin{gather*}
u= \frac{\phi^{2}(r)}{r}\frac{p}{q},
\end{gather*}
and
\begin{itemize}
\itemsep=0pt
\item[(v)]
\begin{gather*}
 \divi{2} (p) \leq \sum\limits_{\ell_{1},\ell_{2} \in
\{0,1\}} 4[\ell_{1}+\ell_{2}\tau]_{2}+2[\ell_{1}+\ell_{2}\tau-h]_{2}
\\
\hphantom{\divi{2} (p) \leq}{} +[\ell_{1}+\ell_{2}\tau+z_{0}-2h]_{2}+[\ell_{1}+\ell_{2}\tau-z_{0}-2h]_{2},
\end{gather*}
\item[(vi)]
\begin{gather*}
 \divi{2} (q)  \leq \sum\limits_{\ell_{1},\ell_{2} \in
\{0,1\}} 2[\ell_{1}+\ell_{2}\tau+h]_{2}+2[\ell_{1}+\ell_{2}\tau+z_{0}+2h]_{2}
\\
\hphantom{\divi{2} (q)  \leq}{}
+2[\ell_{1}+\ell_{2}\tau-z_{0}+2h]_{2}
+[\ell_{1}+\ell_{2}\tau+z_{0}]_{2}+[\ell_{1}+\ell_{2}\tau-z_{0}]_{2},
\end{gather*}
\item[(vii)] $\degr{2} (p)=\degr{2} (q)$,
\item[(viii)] $\weight{2} (p/q)=2\degr{2} (r) h \mod 2\Lambda$.
\end{itemize}
We claim that
\begin{itemize}
\itemsep=0pt
\item[(v$'$)] $\divi{2} (p) \leq \sum\limits_{\ell_{1},\ell_{2} \in \{0,1\}}4[\ell_{1}+\ell_{2}\tau]_{2}$,
\item[(vi$'$)] $\divi{2} (q) \leq \sum\limits_{\ell_{1},\ell_{2} \in
\{0,1\}}[\ell_{1}+\ell_{2}\tau+z_{0}]_{2}+[\ell_{1}+\ell_{2}\tau-z_{0}]_{2}$.
\end{itemize}
Otherwise, arguing as for the proof of the irreducibility of~$G$, we see that (v), (vi) and (viii) would lead to
a~relation of the form
\begin{gather*}
(2\degr{2} (r) + d)h = \ell z_{0} \mod \Lambda
\end{gather*}
for some integer $\ell \in \{-16,\dots,16\}$ and some integer $d>0$ and this would contradict our assumption on $z_{0}$.
Then, (viii) shows that
\begin{gather*}
2\degr{2} (r) h = \ell z_{0} \mod \Lambda
\end{gather*}
for some integer $\ell \in \{-4,\dots,4\}$ and hence $\degr{2} (r)=0$.
So $u=p/q$ with $p,q \in \prodtheta{2}$ satis\-fying~(v$'$) and~(vi$'$) above.
In particular, $-h$ is not a~zero of~$u$.
But $-h$ (which is a~pole of $\phi(a)$) is a~pole of
\begin{gather*}
\phi^{2}(u)+\frac{b}{\phi^{2}(a)}-\phi(a)+\frac{b}{a}.
\end{gather*}
So $-h$ is a~pole of the left hand side of~\eqref{eq10}.
This a~contradiction because $-h$ is not a~pole of the right hand side of~\eqref{eq10}.

Therefore,~$G$ is irreducible and not imprimitive.
So, as explained at the beginning of Section~\ref{sec4}, $G = \{M \in \operatorname{GL}_{2}(\mathbb C) \, \vert \, \det
(M) \in H\}$ where~$H \subset \mathbb C^{\times}$ is the Galois group of $\phi y = b y$, which is easily seen to be
$\mu_{k}$ if~$b$ is a~$k$th root of the unity and $\mathbb C^{\times}$ otherwise.
\end{proof}

\subsection*{Acknowledgements}

Our original interest in dif\/ference equations on elliptic curves arose from discussions with Jean-Pierre Ramis some
years ago.
We thank Jean-Pierre Ramis and Michael Singer for interesting discussions.
We thank the referees for their careful reading and useful suggestions.
The f\/irst author is founded by the labex CIMI.
The second author is partially funded by the French ANR project QDIFF (ANR-2010-JCJC-010501).

\pdfbookmark[1]{References}{ref}

\LastPageEnding

\end{document}